\documentclass[onecolumn, draftcls, 12pt]{IEEEtran}
\usepackage[pdftex]{graphicx}
\usepackage{amsmath}
\usepackage{amssymb}
\usepackage{amsthm}
\usepackage{color}
\usepackage{nccmath}

\makeatletter
\newcommand{\pushright}[1]{\ifmeasuring@#1\else\omit\hfill$\displaystyle#1$\fi\ignorespaces}
\newcommand{\pushleft}[1]{\ifmeasuring@#1\else\omit$\displaystyle#1$\hfill\fi\ignorespaces}
\makeatother

\begin{document}

\title{Stabilization of Networked Control Systems with Sparse Observer-Controller Networks}

\author{Mohammad~Razeghi-Jahromi and Alireza~Seyedi,~\IEEEmembership{Senior Member,~IEEE}
\thanks{M. Razeghi-Jahromi is with the ECE Department, University of Rochester, Rochester,
NY (e-mail: jahromi@ece.rochester.edu).
A. Seyedi is with the EECS Department, University of Central Florida, Orlando,
FL (e-mail: alireza.seyedi@ieee.org).}}

\maketitle
\begin{abstract}
In this paper we provide a set of stability conditions for linear time-invariant networked control systems with arbitrary topology, using a Lyapunov direct approach. We then use these stability conditions to provide a novel low-complexity algorithm for the design of a sparse observer-based control network. We employ distributed observers by employing the output of other nodes to improve the stability of each observer dynamics. To avoid unbounded growth of controller and observer gains, we impose bounds on their norms. The effects of relaxation of these bounds is discussed when trying to find the complete decentralization conditions.
\end{abstract}



\section{Introduction}\label{sec:intro}

Control systems with spatially distributed components have been in use for several decades. In early systems, the components were connected via dedicated hard-wired links carrying the information from the sensors to a central location, where control signals were computed and sent to the actuators. Today, advances in communications technology have enabled us to exchange information via efficient communication networks. These advances have considerably widened the scope of the research on spatially distributed control systems to include communications and network effects explicitly, as they significantly affect the dynamic behavior of the entire system.

\begin{figure}[!t]
\centering
\includegraphics[scale=0.3]{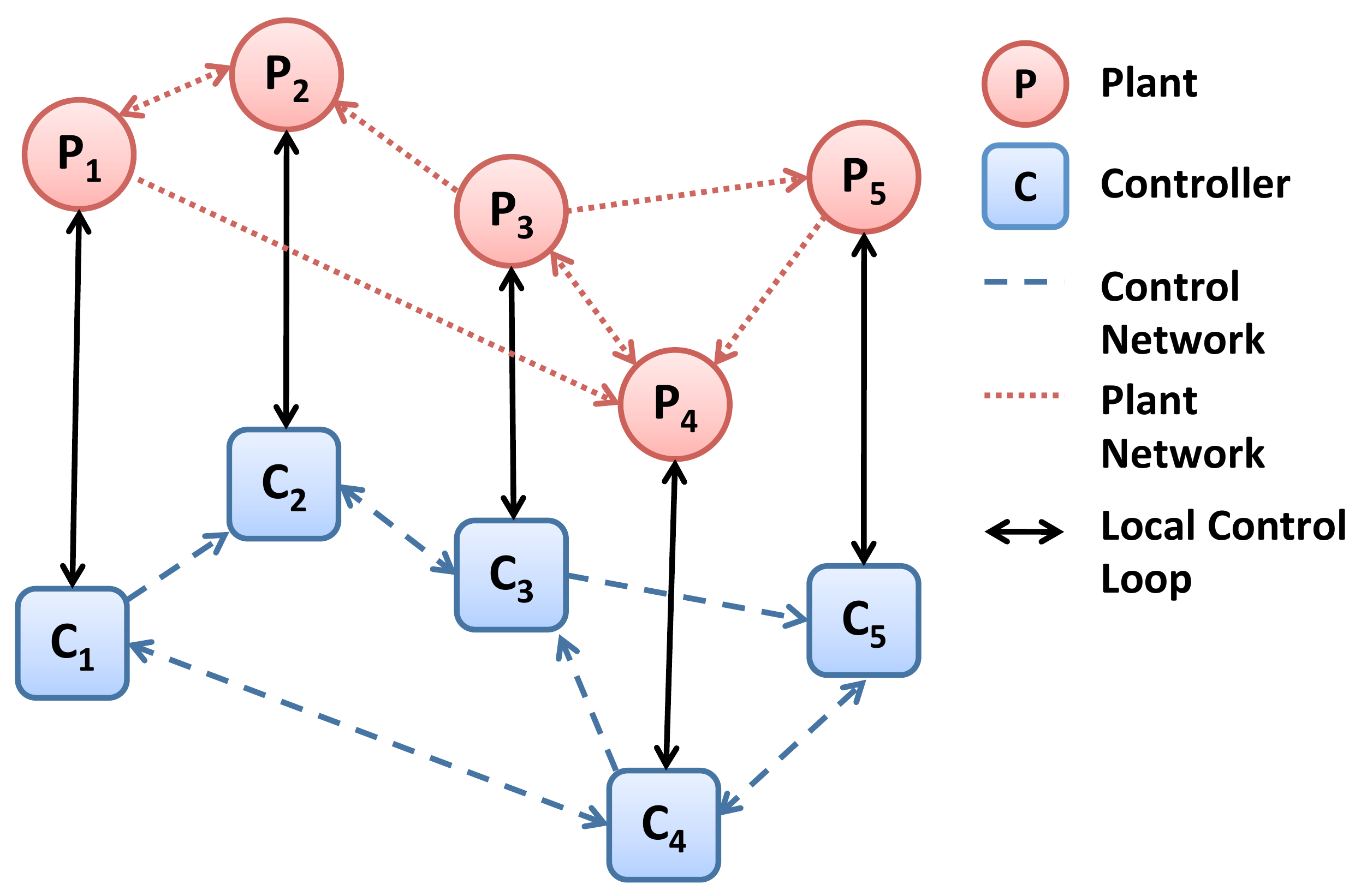}
\vspace{-0.1in}
\caption {A Networked Control System (NCS)}
\vspace{-0.2in}
\label{fig:NCS1}
\end{figure}

Spatially distributed control systems can be abstracted as networked control systems (NCS). An NCS consists of a number of subsystems, each comprising of a plant and a controller, coupled together in a network structure. The interaction of plants with each other forms the {\em plant network}. Control signals are exchanged using the {\em control network}, a.k.a. information, communications, or feedback network (Fig. \ref{fig:NCS1}). Networked control systems have a wide range of applications including electrical networks, transportation networks, factory automation, tele-operations and sensor and actuator networks.

Within this framework, a centralized control approach can be modeled by considering a complete control network, which provides all controllers with access to states of all plants. However, in general, it is not practical to control a large-scale networked system with the centralized approach, where the control law uses the state information of all subsystems, as this requires a large and costly control network. To overcome this limitation, we must resort to decentralized or distributed control strategies \cite{Zhang2001}-\nocite{Baillieul2007,Hespanha2007}\cite{Yue2008}. The decentralized control strategy lies at the opposite end of the spectrum from the centralized approach, where the control law uses only a subsystem's local state information to control the given subsystem. In other words, there is no control network. Such local controls can be effective when the couplings between subsystems are weak \cite{Siljak1991}-\nocite{Siljak2005}\cite{Mazo2010}. However, when the coupling between subsystems are not weak, we may have to use a distributed control approach, where each subsystem uses its own state as well as the state of some other subsystems. This is a middle-of-the-road solution, between centralized and decentralized approaches. Hence, it can achieve stability given stronger subsystem coupling, compared to the decentralized control strategy \cite{Wang2011}, without the complexity of a centralized approach.

Whether a centralized, distributed or decentralized approach is taken, both the dynamics of each subsystem and the network topology, play important roles in the stability of the overall network. It is easy to verify that even if each subsystem is asymptotically stable in isolation, the network may be unstable. In such a scenario to stabilize the NCS a control network carrying state or feedback information between different subsystems may be necessary.


Networked control literature can be classified into two main groups. The first group focuses on the effects of the impairments and limitations imposed by a communication channel, including bandwidth, packet dropout, quantization and delay \cite{montestruque2004}-\nocite{Wang2011,Teel2004,Walsh2002,seiler2005h,yu2005,wu2007}\cite{witrant2007}. The second group, in which this work should be placed, considers the topological and network effects, and investigates how the topology of the plant network affects the overall network behavior, and how a control network can be designed that results in stability or desired performance.

For both decentralized and distributed control approaches, existing works have studied the problem of imposing a priori constraints on communication requirements between subsystems.
It has been shown that under a structural condition, namely {\em quadratic invariance}, finding optimal controllers can be cast as a convex optimization problem \cite{Lall2006}-\nocite{Rotkowitz2009,Rotkowitz2010}\cite{Swigart2009}.  Other work have shown similar results, conditioned on the network being a partially ordered set (poset) \cite{Shah2011,Parrilo2011}. This constraint is closely related to quadratic invariance, however, it can lead to more computationally efficient solutions and provides better insight into the topology of the optimal controllers.

While these results are elegant and important, they impose restrictions on the topology of the plant network. For networks with arbitrary topology, the key question concerning the design of the control network is one of {\em topological information requirements} and can be framed as: {\em Which nodes should be given the state and output information of a particular node, in order for the local controllers to be able to satisfy a global control objective?} This question is critical in the design of massively distributed control systems, such as the Smart Grid \cite{Massoud2011}-\nocite{Schuppen2011,Lafortune2003}\cite{Lin2007}.

In addressing this key question, the goal is often to find the sparsest control network that satisfies the requirements. This problem has been considered in various settings and solutions have been proposed that can be used to find suboptimal controllers. In \cite{lin2011,fardad2011sparsity}, a non-convex condition is proposed which may be solved numerically. It should be mentioned that different notions of sparsity have been employed. For instance, in \cite{fardad2011sparsity} the total number of non-zero elements in the coupling matrices is considered as a measure of sparsity, as opposed to the number of non-zero coupling matrices (number of links in the network).

In this paper, we first develop a set of stability conditions that guarantee global asymptotic stability, using the Lyapunov direct method. These conditions are significantly less conservative than our prior results in \cite{razeghi2011} and include state estimation, among other improvements. We then use these conditions to explore the problem of designing a sparse observer-controller network for a given plant network with arbitrary topology. We take a broader look at the topological information requirements by taking into account the distributed state estimation problem, generally neglected by the existing works. We proceed to provide a solution for finding the sparsest observer-controller network that satisfies our set of stability conditions. We show that stabilization of NCS that are partially ordered set is trivial under our condition.

We consider a linear time-invariant (LTI) NCS with arbitrary topology and provide a methodology to design a sparse observer and controller network. We assume that communication links do not have any bandwidth limitation, data loss or induced network delays. Moreover, we find the conditions and the bounds on norm of local controller gain matrices which make complete decentralization is possible.

\section{Notation and Problem Definition}\label{sec:notation}

\subsubsection{Notation}
Matrices and vectors are denoted by capital and lower-case bold letters, respectively. We use $\mathbb{S}^{n}_{++}$ to denote the set of real symmetric positive definite $n\times n$ matrices. Generalized matrix inequality, $\prec$, is defined by the positive definite cone between symmetric matrices. The Euclidean ($l_2$) vector norm and the induced $l_2$ matrix norm are represented by $\| \cdot \|$. By $\lambda_{\mbox{\footnotesize min}}(.)$, $\lambda_{\mbox{\footnotesize max}}(.)$ and $\sigma_{\mbox{\footnotesize max}}(.)$ we denote the smallest and largest eigenvalue and the largest singular value of the argument, respectively. The $m\times n$ unit matrix consisting of all ones is denoted by $\mathbf{1}^{m\times n}$. We let $\mathcal{N}=\{1,\ldots,N\}$ and $\mathcal{N}_i = \mathcal{N}-\{i\}$. The indicator function of $x$ is represented by $1_x$ and column-stacking operator is denoted by $\mbox{vec}(.)$.

\subsubsection{Problem Definition}
Consider a network of $N$ coupled LTI subsystems. The state of the $i$th plant, $\mathbf{x}_i(t)$, is given by
\begin{eqnarray}
\label{eq: 1}
\dot{\mathbf{x}}_i(t) & = & \mathbf{A}_i\mathbf{x}_i(t)+\mathbf{B}_i\mathbf{u}_i(t)+\sum_{j\in \mathcal{N}_i}\mathbf{H}_{ij}\mathbf{x}_j(t)\nonumber\\
\mathbf{y}_i(t) & = & \mathbf{C}_i\mathbf{x}_i(t),
\end{eqnarray}
where $\mathbf{u}_i(t)$ and $\mathbf{y}_i(t)$ are input and output of the $i$th subsystem, and $\mathbf{A}_i$, $\mathbf{B}_i$, $\mathbf{C}_i$ and $\mathbf{H}_{ij}$ are known matrices. We assume that $(\mathbf{A}_i,\mathbf{B}_i)$ are controllable and $(\mathbf{A}_i,\mathbf{C}_i)$ are observable. We consider an arbitrary directed network without self-loops. That is, $\mathbf{H}_{ii}=\mathbf{0}$, and $\mathbf{H}_{ij}$ and $\mathbf{H}_{ji}$ are not necessarily equal. We look for a distributed stabilizing observer-based controller of the form
\begin{eqnarray}
\label{eq: 2}
\dot{\hat{\mathbf{x}}}_i(t) & = & \mathbf{A}_i\hat{\mathbf{x}}_i(t)+\mathbf{B}_i\mathbf{u}_i(t)+\sum_{j\in \mathcal{N}_i}\mathbf{H}_{ij}\hat{\mathbf{x}}_j(t)+\mathbf{M}_i(\mathbf{C}_i\hat{\mathbf{x}}_i(t)-\mathbf{y}_i(t))+
\sum_{j\in \mathcal{N}_i}\mathbf{O}_{ij}(\mathbf{C}_j\hat{\mathbf{x}}_j(t)-\mathbf{y}_j(t)),\nonumber\\
\mathbf{u}_i(t) & = & \mathbf{K}_i\hat{\mathbf{x}}_i(t)+\sum_{j\in \mathcal{N}_i}\mathbf{L}_{ij}\hat{\mathbf{x}}_j(t),
\end{eqnarray}
where $\hat{\mathbf{x}}_i(t)$ is the estimate of $\mathbf{x}_i(t)$, $\mathbf{K}_i$ and $\mathbf{L}_{ij}$ are local and coupling controller gains, and $\mathbf{M}_i$ and $\mathbf{O}_{ij}$ are local and coupling observer gains, respectively. Note that to estimate $\mathbf{x}_i(t)$, we not only use output of subsystem $i$, but also outputs of (potentially) all other subsystems. This is dual to the concept of distributed control. Our objective is to find distributed observer-based control law (\ref{eq: 2}), using feedback from (potentially) all other subsystems to stabilize the plant network with a sparse control network. That is, we aim to find $\mathbf{K}_i,\mathbf{M}_i,\mathbf{L}_{ij}$ and $\mathbf{O}_{ij}$, such that the overall network is globally asymptotically stable and that the number of links in the control network (number of non-zero coupling gains $\mathbf{L}_{ij}$ and $\mathbf{O}_{ij}$) is minimized. We also impose constraints
\begin{eqnarray}
\label{eq: 5}
\| \mathbf{K}_i\| \leq \kappa_i,~\| \mathbf{M}_i\| \leq \mu_i,~\| \mathbf{L}_{ij}\|\leq \iota_{ij},~\| \mathbf{O}_{ij}\|\leq \omega_{ij},
\end{eqnarray}
to avoid undesirably large gains.

Defining $\mathbf{x}(t)=\mbox{vec}(\mathbf{x}_i(t))$, $\mathbf{u}(t)=\mbox{vec}(\mathbf{u}_i(t))$, $\mathbf{y}(t)=\mbox{vec}(\mathbf{y}_i(t))$, (\ref{eq: 1}) reduces to
\begin{eqnarray}
\label{eq: 1-1}
\dot{\mathbf{x}}(t) = \mathbf{A}\mathbf{x}(t)+\mathbf{B}\mathbf{u}(t)+\mathbf{H}\mathbf{x}(t),\,\,\,\,\mbox{and}\,\,\,\,\mathbf{y}(t) = \mathbf{C}\mathbf{x}(t),
\end{eqnarray}
where $\mathbf{A}=\mbox{diag}(\mathbf{A}_i)$, $\mathbf{B}=\mbox{diag}(\mathbf{B}_i)$, $\mathbf{C}=\mbox{diag}(\mathbf{C}_i)$ and $\mathbf{H}=\left[\mathbf{H}_{ij}\right]$.
Moreover, (\ref{eq: 2}) yields
\begin{eqnarray}
\label{eq: 1-3}
\dot{\hat{\mathbf{x}}}(t) & = & \mathbf{A}\hat{\mathbf{x}}(t)+\mathbf{B}\mathbf{u}(t)+\mathbf{H}\hat{\mathbf{x}}(t)
+\mathbf{M}(\mathbf{C}\hat{\mathbf{x}}(t)-\mathbf{y}(t))+\mathbf{O}(\mathbf{C}\hat{\mathbf{x}}(t)-\mathbf{y}(t)),\nonumber\\
\mathbf{u}(t) & = & \mathbf{K}\hat{\mathbf{x}}(t)+\mathbf{L}\hat{\mathbf{x}}(t),
\end{eqnarray}
where $\mathbf{K}=\mbox{diag}(\mathbf{K}_i)$, $\mathbf{M}=\mbox{diag}(\mathbf{M}_i)$, $\mathbf{L}=\left[\mathbf{L}_{ij}\right]$, with $\mathbf{L}_{ii}=\mathbf{0}$ and $\mathbf{O}=\left[\mathbf{O}_{ij}\right]$, with $\mathbf{O}_{ii}=\mathbf{0}$.\

Defining error $\mathbf{e}(t)\triangleq \hat{\mathbf{x}}(t)-\mathbf{x}(t)$ reduces (\ref{eq: 1-1}) and (\ref{eq: 1-3}) to
\begin{align}
\label{eq: 3}
&\dot{\mathbf{x}}(t)=\left[\mathbf{A}+\mathbf{H}+\mathbf{B}(\mathbf{K}+\mathbf{L})\right]\mathbf{x}(t)+\mathbf{B}(\mathbf{K}+\mathbf{L})\mathbf{e}(t),\\
\label{eq: 4}
&\dot{\mathbf{e}}(t)=\left[\mathbf{A}+\mathbf{H}+(\mathbf{M}+\mathbf{O})\mathbf{C}\right]\mathbf{e}(t).
\end{align}
This is a networked linear cascade dynamical system with the equilibrium point $(\mathbf{x},\mathbf{e})=(\mathbf{0},\mathbf{0})$.

\section{Network Stability Conditions}\label{sec:stability}

The following lemma provides conditions for the globally asymptotic stability of the network.

\emph{Lemma 1}: The equilibrium point of the system in (\ref{eq: 3}) and (\ref{eq: 4}), which is $(\mathbf{x},\mathbf{e})=(\mathbf{0},\mathbf{0})$, is globally asymptotically stable, if there exist $\mathbf{K}$, $\mathbf{L}$, $\mathbf{P}$, $\mathbf{M}$, $\mathbf{O}$ and $\hat{\mathbf{P}}$ such that
\begin{eqnarray}
\begin{array}{rl}
\label{eq: 7}
\left[\mathbf{A}+\mathbf{H}+\mathbf{B}(\mathbf{K}+\mathbf{L})\right]^{T}\mathbf{P}+\mathbf{P}\left[\mathbf{A}+\mathbf{H}+\mathbf{B}(\mathbf{K}+\mathbf{L})\right]
+2\boldsymbol\beta\circ \mathbf{P}\prec \mathbf{0} & (\mbox{S1})\\
\left[\mathbf{A}+\mathbf{H}+(\mathbf{M}+\mathbf{O})\mathbf{C}\right]^{T}\hat{\mathbf{P}}+\hat{\mathbf{P}}\left[\mathbf{A}+\mathbf{H}+(\mathbf{M}+\mathbf{O})\mathbf{C}\right]+2\boldsymbol\beta\circ \hat{\mathbf{P}}\prec \mathbf{0} & (\mbox{S2})\\
\mathbf{P}=\mbox{diag}(\mathbf{P}_i) \succ \mathbf{0} & (\mbox{S3})\\
\hat{\mathbf{P}}=\mbox{diag}(\hat{\mathbf{P}}_i)\succ  \mathbf{0} & (\mbox{S4})
\end{array}
\end{eqnarray}
where $\boldsymbol\beta=\mbox{diag}\left(\beta_i\mathbf{1}^{n_i\times n_i}\right)$, $\beta_i\geq 0$ is the stability margin for subsystem $i$, and $\circ$ is the  Hadamard product.

\begin{IEEEproof}
(S1) and (S3) guarantee input-to-state stability of (\ref{eq: 3}) while (S2) and (S4) guarantee global asymptotical stability of (\ref{eq: 4}). Since the system is a cascaded dynamical system, $(\mathbf{x},\mathbf{e})=(\mathbf{0},\mathbf{0})$, is globally asymptotically stable \cite{haddad2011}.
\end{IEEEproof}

\section{Sparse Control Network Design}\label{sec:design}
Now, our objective is to design a control network with minimum number of links that satisfies stability condition (\ref{eq: 7}), under gain constraints (\ref{eq: 5}). This problem can be formulated as
\begin{eqnarray}
\label{eq: 39.11}
\begin{array}{cl}
\underset{\mathbf{K},\mathbf{L},\mathbf{P},\mathbf{M},\mathbf{O},\hat{\mathbf{P}}}{\text{minimize}}&\sum_{i,j\in \mathcal{N}}1_{\left\{\mathbf{L}_{ij}\neq \mathbf{0}~\mbox{\footnotesize{or}}~\mathbf{O}_{ij}\neq \mathbf{0}\right\}}\\
\mbox{subject to} & \mbox{(\ref{eq: 5}) and (\ref{eq: 7})}.
\end{array}
\end{eqnarray}

Unfortunately, besides the fact that the objective function is integer valued, the first two constraints in (\ref{eq: 7}) are not convex. In the following we convexify this problem by restricting its domain.

\emph{Theorem 1}: System (\ref{eq: 1}) with controller (\ref{eq: 2}) is globally asymptotically stable, and bounds (\ref{eq: 5}) are satisfied, if the following convex mixed-binary program has a solution
\begin{eqnarray}
\label{eq: 39}
\begin{array}{crl}
\mbox{minimize} & \sum_{i,j\in \mathcal{N}}\alpha_{ij} \quad\quad\quad\quad\quad\quad\quad\quad\quad\quad\quad\quad&\\
\mbox{subject to}& \mathbf{F}(\mathbf{Z},\mathbf{W},\mathbf{Y},\boldsymbol{\alpha})+\mathbf{F}^{T}(\mathbf{Z},\mathbf{W},\mathbf{Y},\boldsymbol{\alpha})\prec \mathbf{0} & \mbox{(C1)}\\
&\hat{\mathbf{F}}(\hat{\mathbf{P}},\hat{\mathbf{W}},\hat{\mathbf{Y}},\hat{\boldsymbol{\alpha}})+\hat{\mathbf{F}}^{T}(\hat{\mathbf{P}},\hat{\mathbf{W}},\hat{\mathbf{Y}},\hat{\boldsymbol{\alpha}})\prec \mathbf{0} & \mbox{(C2)}\\
&\mathbf{Z} \succ \mathbf{0} & \mbox{(C3)}\\
&\hat{\mathbf{P}} \succ \mathbf{0} & \mbox{(C4)}\\
&\kappa_i\lambda_{\mbox{\footnotesize min}}(\mathbf{Z}_i)-\sigma_{\mbox{\footnotesize max}}(\mathbf{W}_i) \geq  0 & \mbox{(C5)}\\
&\iota_{ij}\lambda_{\mbox{\footnotesize min}}(\mathbf{Z}_j)-\sigma_{\mbox{\footnotesize max}}(\mathbf{Y}_{ij}) \geq  0 & \mbox{(C6)}\\
&\mu_i\lambda_{\mbox{\footnotesize min}}(\hat{\mathbf{P}}_i)-\sigma_{\mbox{\footnotesize max}}(\hat{\mathbf{W}}_i) \geq 0 & \mbox{(C7)}\\
&\omega_{ij}\lambda_{\mbox{\footnotesize min}}(\hat{\mathbf{P}}_i)-\sigma_{\mbox{\footnotesize max}}(\hat{\mathbf{Y}}_{ij}) \geq  0 & \mbox{(C8)}\\
&\alpha_{ik}\in \{0,1\} & \mbox{(C9)}
\end{array}
\end{eqnarray}
for all $i,k\in \mathcal{N}$, where
\begin{eqnarray}
\mathbf{F}(\mathbf{Z},\mathbf{W},\mathbf{Y},\boldsymbol{\alpha})=(\mathbf{A}+\mathbf{H})\mathbf{Z}+\mathbf{B}(\mathbf{W}+\boldsymbol\alpha\circ\mathbf{Y})+\boldsymbol \beta\circ \mathbf{Z},\nonumber\\
\hat{\mathbf{F}}(\hat{\mathbf{P}},\hat{\mathbf{W}},\hat{\mathbf{Y}},\hat{\boldsymbol{\alpha}})=\hat{\mathbf{P}}(\mathbf{A}+\mathbf{H})+(\hat{\mathbf{W}}
+\hat{\boldsymbol\alpha}\circ\hat{\mathbf{Y}})\mathbf{C}+\boldsymbol \beta\circ \hat{\mathbf{P}},\nonumber
\end{eqnarray}
where $\mathbf{Z}=\mbox{diag}(\mathbf{Z}_i)$, $\hat{\mathbf{P}}=\mbox{diag}(\hat{\mathbf{P}}_i)$, $\mathbf{W}=\mbox{diag}(\mathbf{W}_i)$, $\hat{\mathbf{W}}=\mbox{diag}(\hat{\mathbf{W}}_i)$, $\mathbf{Y}=[\mathbf{Y}_{ij}]$, $\hat{\mathbf{Y}}=[\hat{\mathbf{Y}}_{ij}]$ with $\mathbf{Y}_{ii}=\hat{\mathbf{Y}}_{ii}=\mathbf{0}$, $\boldsymbol\alpha=\left[\alpha_{ij}\mathbf{1}_{ij}^{m_i\times n_j}\right]$, $\hat{\boldsymbol\alpha}=\left[\alpha_{ij}\mathbf{1}_{ij}^{n_i\times r_j}\right]$ and $\alpha_{ii}=0$. Furthermore, if $(\mathbf{Z}^\star_i,\hat{\mathbf{P}}^\star_i,\mathbf{W}^\star_i,\hat{\mathbf{W}}^\star_i,\mathbf{Y}^\star_{ij},\hat{\mathbf{Y}}^\star_{ij},\alpha^\star_{ij})$ is a solution of \eqref{eq: 39}, the controller and observer gains are
\begin{eqnarray}
\label{eq: 39.1}
\begin{array}{lll}
\mathbf{K}_i=\mathbf{W}^\star_i(\mathbf{Z}^\star_i)^{-1}, & \mathbf{L}_{ij}=\alpha^\star_{ij}\mathbf{L}^\star_{ij},  & \mathbf{L}^\star_{ij}=\mathbf{Y}^{\star}_{ij}(\mathbf{Z}^\star_j)^{-1}\\
\mathbf{M}_i=(\hat{\mathbf{P}}^\star_i)^{-1}\hat{\mathbf{W}}^\star_i & \mathbf{O}_{ij}=\alpha^\star_{ij}\mathbf{O}^\star_{ij}, & \mathbf{O}^\star_{ij}=(\hat{\mathbf{P}}^\star_i)^{-1}\hat{\mathbf{Y}}^{\star}_{ij}.
\end{array}
\end{eqnarray}

\begin{IEEEproof}
By defining $\mathbf{Z}\triangleq  \mathbf{P}^{-1}$, $\mathbf{W}\triangleq  \mathbf{K}\mathbf{P}^{-1}$, and $\mathbf{Y}\triangleq \mathbf{L}\mathbf{P}^{-1}$, we can write (S1) and (S3) as
\begin{align}
\label{eq: 19}
\left[(\mathbf{A}+\mathbf{H})\mathbf{Z}+\mathbf{B}(\mathbf{W}+\mathbf{Y})\right]^T+\left[(\mathbf{A}+\mathbf{H})\mathbf{Z}+\mathbf{B}(\mathbf{W}+\mathbf{Y})\right]+2\boldsymbol \beta \circ \mathbf{Z}&\prec \mathbf{0},\nonumber\\
\mathbf{Z}=\mbox{diag}(\mathbf{Z}_i)&\succ \mathbf{0}.
\end{align}
The original variables can then be found from $\mathbf{P}=\mathbf{Z}^{-1}$, $\mathbf{K}=\mathbf{W}\mathbf{Z}^{-1}$ and $\mathbf{L}=\mathbf{Y}\mathbf{Z}^{-1}$.

To design a sparse control network, we seek a set of $\mathbf{L}_{ij}$ that guarantee stability, with minimum number of non-zero $\mathbf{L}_{ij}$. Note that when the $ij$th link is included in the control network, we should use the gain $\mathbf{L}^{\star}_{ij}$. On the other hand when the $ij$th link is not used we have $\mathbf{L}_{ij}=\mathbf{0}$. This can be expressed as $\mathbf{L}_{ij}=\alpha_{ij}\mathbf{L}^{\star}_{ij}$ or for overall system $\mathbf{L}=\boldsymbol\alpha\circ\mathbf{L}^{\star}$, where $\alpha_{ij}=1$ if the $ij$th link is used in the control network and $\alpha_{ij}=0$ if it is not.

To find the variables $\mathbf{Z}$, $\mathbf{W}$ and $\mathbf{Y}$ in (\ref{eq: 19}), we can rewrite (\ref{eq: 19}) as $\alpha_{ik}\in \{0,1\}$, $\mathbf{Z}\succ \mathbf{0}$ and $\mathbf{F}(\mathbf{Z},\mathbf{W},\mathbf{Y},\boldsymbol{\alpha})+\mathbf{F}^{T}(\mathbf{Z},\mathbf{W},\mathbf{Y},\boldsymbol{\alpha})\prec \mathbf{0},$ which are (C1), (C3), and (C9) in (\ref{eq: 39}).

While changing the variables from $\mathbf{P}$, $\mathbf{K}$ and $\mathbf{L}$ to $\mathbf{Z}$, $\mathbf{W}$ and $\mathbf{Y}$ convexified the first two constraints in (\ref{eq: 7}), it caused the local gain constraint (\ref{eq: 5}) to become non-convex. To remedy this, we can convexify this constraint by first upper bounding the norm of $\mathbf{K}_i$ as
\begin{eqnarray}
\label{eq: 27}
\| \mathbf{K}_i\|  =  \| \mathbf{W}_i\mathbf{Z}^{-1}_i\| \leq \| \mathbf{W}_i\|\| \mathbf{Z}^{-1}_i\|  = \sigma_{\mbox{\footnotesize max}}(\mathbf{W}_i)\lambda_{\mbox{\footnotesize max}}(\mathbf{Z}^{-1}_i) = \frac{\sigma_{\mbox{\footnotesize max}}(\mathbf{W}_i)}{\lambda_{\mbox{\footnotesize min}}(\mathbf{Z}_i)},
\end{eqnarray}
and forcing (\ref{eq: 5}) by upper bounding (\ref{eq: 27}) by $\kappa_i$:
\begin{align}
\label{eq: 28}
\| \mathbf{K}_i\|\leq \frac{\sigma_{\mbox{\footnotesize max}}(\mathbf{W}_i)}{\lambda_{\mbox{\footnotesize min}}(\mathbf{Z}_i)}\leq \kappa_i.
\end{align}
Equivalently $\kappa_i\lambda_{\mbox{\footnotesize min}}(\mathbf{Z}_i)-\sigma_{\mbox{\footnotesize max}}(\mathbf{W}_i)\geq 0$,
which is a convex constraint. Similarly, we have $\iota_{ij}\lambda_{\mbox{\footnotesize min}}(\mathbf{Z}_j)-\sigma_{\mbox{\footnotesize max}}(\mathbf{Y}_{ij})\geq 0$. These provide constraints (C5) and (C6) in (\ref{eq: 39}).\

Similarly, we can convexify the general stability condition (S2) with given bounds in (\ref{eq: 5}) as constraints (C2), (C4), (C7) and (C8) in (\ref{eq: 39}) where $\hat{\mathbf{W}}\triangleq  \hat{\mathbf{P}}\mathbf{M}$ and $\hat{\mathbf{Y}}\triangleq \hat{\mathbf{P}} \mathbf{O}$. Note that if link $ij$ is used it can carry both $\mathbf{L}_{ij}\hat{\mathbf{x}}_j(t)$ and $\mathbf{O}_{ij}\mathbf{y}_j(t)$. Thus the same variables $\alpha_{ij}$ should be used for both controller and observer links. Consequently, $\mathbf{O}_{ij}=\alpha_{ij}\mathbf{O}^{\star}_{ij}$ or $\mathbf{O}=\hat{\boldsymbol\alpha}\circ\mathbf{O}^{\star}$.

Minimizing the number of communication links is equivalent to minimizing the number of $\alpha_{ij}=1$, or in other words, minimizing the sum of $\alpha_{ij}$ subject to constraints in (\ref{eq: 39}).
\end{IEEEproof}

\emph{Theorem 2}: The controlled network (\ref{eq: 3}) and (\ref{eq: 4}) are stable if there are no constraints on the norm of gain matrices and
\begin{eqnarray}
\label{eq: 253}
\|\mathbf{B}\mathbf{L}+\mathbf{H}\|<\frac{\lambda_{\mbox{\footnotesize min}}(\mathbf{Q})}{2\lambda_{\mbox{\footnotesize max}}(\mathbf{P})}-\beta_{\max},\\
\label{eq: 253-1}
\|\mathbf{O}\mathbf{C}+\mathbf{H}\|<\frac{\lambda_{\mbox{\scriptsize min}}(\hat{\mathbf{Q}})}{2\lambda_{\mbox{\scriptsize max}}(\hat{\mathbf{Z}})}-\beta_{\max},
\end{eqnarray}
where $\beta_{\max}\triangleq \underset{i}{\max}~\beta_i$ and positive definite matrices $\mathbf{P}$, $\mathbf{Q}=\mbox{diag}(\mathbf{Q}_i)$, $\hat{\mathbf{Z}}$ and $\hat{\mathbf{Q}}=\mbox{diag}(\hat{\mathbf{Q}}_i)$ are the solution of Algebraic Riccati Equations (AREs)
\begin{eqnarray}
\label{eq: 29-2}
\mathbf{A}^T\mathbf{P}+\mathbf{P}\mathbf{A}-\mathbf{P}\mathbf{B}\mathbf{B}^T\mathbf{P}+\mathbf{Q}&=&\mathbf{0},\\
\label{eq: 29-3}
\hat{\mathbf{Z}}\mathbf{A}^T+\mathbf{A}\hat{\mathbf{Z}}-\hat{\mathbf{Z}}\mathbf{C}^T\mathbf{C}\hat{\mathbf{Z}}+\hat{\mathbf{Q}}&=&\mathbf{0}.
\end{eqnarray}

\begin{IEEEproof}
Consider feedback $\mathbf{K}=-\frac{1}{2}\mathbf{B}^T\mathbf{P}$ in (\ref{eq: 7}), we have
\begin{align}
\label{eq: 349-1}
(\mathbf{A}+\mathbf{H}+\mathbf{B}\mathbf{L})^T\mathbf{P}+\mathbf{P}(\mathbf{A}+\mathbf{H}+\mathbf{B}\mathbf{L})-\mathbf{P}\mathbf{B}\mathbf{B}^T\mathbf{P}+2\boldsymbol \beta \circ \mathbf{P}&\prec \mathbf{0},\nonumber\\
\mathbf{P}=\mbox{diag}(\mathbf{P}_i)&\succ \mathbf{0}.
\end{align}
Since $(\mathbf{A},\mathbf{B})$ is controllable, there exist $\mathbf{P},\mathbf{Q} \in \mathbb{S}^{n}_{++}$ such that
\begin{align}
\label{eq: 354-1}
\mathbf{A}^T\mathbf{P}+\mathbf{P}\mathbf{A}-\mathbf{P}\mathbf{B}\mathbf{B}^T\mathbf{P}+\mathbf{Q}=\mathbf{0}.
\end{align}
Substituting $\mathbf{Q}$ from (\ref{eq: 354-1}) to (\ref{eq: 349-1}), we need to have
\begin{align}
\label{eq: 355-1}
(\mathbf{B}\mathbf{L}+\mathbf{H})^T\mathbf{P}+\mathbf{P}(\mathbf{B}\mathbf{L}+\mathbf{H})+2\boldsymbol \beta \circ \mathbf{P}-\mathbf{Q}\prec\mathbf{0},
\end{align}
To satisfy (\ref{eq: 355-1}), it is sufficient to have
\begin{align}
\lambda_{\mbox{\footnotesize max}}\left[(\mathbf{B}\mathbf{L}+\mathbf{H})^T\mathbf{P}+\mathbf{P}(\mathbf{B}\mathbf{L}+\mathbf{H})+2\boldsymbol \beta \circ \mathbf{P}\right]<\lambda_{\mbox{\footnotesize min}}(\mathbf{Q}).\nonumber
\end{align}
Since $| \lambda_i(\mathbf{G})|\leq \|\mathbf{G}\|$, we can upper bound the left hand side
\begin{align}
\label{eq: 351-2}
\lambda_{\mbox{\footnotesize max}}\left[(\mathbf{B}\mathbf{L}+\mathbf{H})^T\mathbf{P}+\mathbf{P}(\mathbf{B}\mathbf{L}+\mathbf{H})+2\boldsymbol \beta \circ \mathbf{P}\right]&\leq \|(\mathbf{B}\mathbf{L}+\mathbf{H})^T\mathbf{P}+\mathbf{P}(\mathbf{B}\mathbf{L}+\mathbf{H})+2\boldsymbol \beta \circ \mathbf{P}\|\nonumber\\
&\leq 2\|\mathbf{P}\|\left(\|\mathbf{B}\mathbf{L}+\mathbf{H}\|+\ \beta_{\max}\right),
\end{align}
Thus, to satisfy (\ref{eq: 355-1}), it suffices to have
\begin{align}
\label{eq: 352-1}
\|\mathbf{B}\mathbf{L}+\mathbf{H}\|<\frac{\lambda_{\mbox{\footnotesize min}}(\mathbf{Q})}{2\lambda_{\mbox{\footnotesize max}}(\mathbf{P})}-\beta_{\max}.
\end{align}
Similarly norm bound in (\ref{eq: 253-1}) is obtained by setting $\mathbf{M}=-\frac{1}{2}\hat{\mathbf{P}}^{-1}\mathbf{C}^T$ and $\hat{\mathbf{Z}}\triangleq\hat{\mathbf{P}}^{-1}$ in (\ref{eq: 7}).
\end{IEEEproof}

Compared to Theorem 1, Theorem 2 is more conservative. However, it provides some insight into decentralized control as described in the following.

\emph{Theorem 3}: Decentralized control is possible if $\kappa_i\geq\underline{\kappa_i}$ and $\mu_i\geq\underline{\mu_i}$ and we have either
\begin{itemize}
\item{$\mathbf{B}\mathbf{B}^T$ is non-singular, or}
\item{$\|\mathbf{H}\|<\frac{\lambda_{\mbox{\scriptsize min}}(\mathbf{Q})}{2\lambda_{\mbox{\scriptsize max}}(\mathbf{P})}-\beta_{\max}$, where $\mathbf{P}\in \mathbb{S}^{n}_{++}$ and $\mathbf{Q}\in \mathbb{S}^{n}_{++}$ are the solution of (\ref{eq: 29-2})
}
\end{itemize}
and, either
\begin{itemize}
\item{$\mathbf{C}^T\mathbf{C}$ is non-singular, or}
\item{$\|\mathbf{H}\|<\frac{\lambda_{\mbox{\scriptsize min}}(\hat{\mathbf{Q}})}{2\lambda_{\mbox{\scriptsize max}}(\hat{\mathbf{Z}})}-\beta_{\max}$, where $\hat{\mathbf{Z}}\in \mathbb{S}^{n}_{++}$ and $\hat{\mathbf{Q}}\in \mathbb{S}^{n}_{++}$ are the solution of (\ref{eq: 29-3}).
}
\end{itemize}
The bounds on decentralized controller gains are
\begin{align}
\label{eq: 239}
\underline{\kappa_i}= \frac{1}{2}\|\mathbf{B}_i^T\mathbf{Z}_i^{-1}\|,~~  \underline{\mu_i}= \frac{1}{2}\|\hat{\mathbf{P}}_i^{-1}\mathbf{C}_i^T\|,
\end{align}
where $\mathbf{Z}_i$ and $\hat{\mathbf{P}}_i$ are the solution of

\begin{eqnarray}
\label{eq: 39-1}
\begin{array}{crc}
\text{maximize}&\sum_{i\in \mathcal{N}}\left(\lambda_{\mbox{\footnotesize min}}(\mathbf{Z}_i)+\lambda_{\mbox{\footnotesize min}}(\hat{\mathbf{P}}_i)\right)\quad\quad\quad\quad\quad\quad\quad\quad&\\
\mbox{subject to}&\mathbf{Z}(\mathbf{A}+\mathbf{H})^T+(\mathbf{A}+\mathbf{H})\mathbf{Z}+2\boldsymbol \beta \circ \mathbf{Z}-\mathbf{B}\mathbf{B}^T\prec \mathbf{0}&\mbox{(D1)}\\
&(\mathbf{A}+\mathbf{H})^T\hat{\mathbf{P}}+\hat{\mathbf{P}}(\mathbf{A}+\mathbf{H})+2\boldsymbol \beta \circ \hat{\mathbf{P}}-\mathbf{C}^T\mathbf{C}\prec  \mathbf{0} & \mbox{(D2)}\\
&\mathbf{Z}=\mbox{diag}(\mathbf{Z}_i) \succ \mathbf{0} & \mbox{(D3)}\\
&\hat{\mathbf{P}}=\mbox{diag}(\hat{\mathbf{P}}_i) \succ  \mathbf{0} & \mbox{(D4)}
\end{array}
\end{eqnarray}
which is a convex and always feasible under the above conditions.

\begin{IEEEproof}
Set $\mathbf{L}=\mathbf{0}$. Pre- and post-multiplication of (\ref{eq: 349-1}) by $\mathbf{P}^{-1}=\mathbf{Z}$ yields
\begin{align}
\label{eq: 350}
\mathbf{Z}(\mathbf{A}+\mathbf{H})^T+(\mathbf{A}+\mathbf{H})\mathbf{Z}+2\boldsymbol \beta \circ \mathbf{Z}-\mathbf{B}\mathbf{B}^T&\prec \mathbf{0},\nonumber\\
\mathbf{Z}=\mbox{diag}(\mathbf{Z}_i)&\succ \mathbf{0}.
\end{align}
To satisfy (\ref{eq: 350}), it is sufficient to have
\begin{align}
\label{eq: 352}
\lambda_{\mbox{\footnotesize max}}(\mathbf{Z})&<\frac{\frac{1}{2}\lambda_{\mbox{\footnotesize min}}(\mathbf{\mathbf{B}\mathbf{B}^T})}{\|\mathbf{A}+\mathbf{H}\|+\beta_{\max}},\nonumber\\
\mathbf{Z}=\mbox{diag}(\mathbf{Z}_i)&\succ \mathbf{0}.
\end{align}
If $\mathbf{B}\mathbf{B}^T$ is non-singular, (\ref{eq: 352}) can be satisfied. If $\mathbf{B}\mathbf{B}^T$ is singular, we have $\lambda_{\mbox{\footnotesize min}}(\mathbf{\mathbf{B}\mathbf{B}^T})=0$. Using Theorem 2, we need to have
\begin{align}
\label{eq: 356}
\|\mathbf{H}\|<\frac{\lambda_{\mbox{\footnotesize min}}(\mathbf{Q})}{2\lambda_{\mbox{\footnotesize max}}(\mathbf{P})}-\beta_{\max},
\end{align}
where $\mathbf{P}$ and $\mathbf{Q}$ are the solution of (\ref{eq: 29-2}). Therefore, if either $\mathbf{B}\mathbf{B}^T$ is non-singular or (\ref{eq: 356}) is satisfied, (D1) and (D3) in (\ref{eq: 39-1}) are always feasible. Moreover, we have
\begin{eqnarray}
\label{eq: 251}
\min~\sum_{i\in \mathcal{N}}\|\mathbf{K}_i\| =\min~\frac{1}{2}\sum_{i\in \mathcal{N}}\left\|\mathbf{B}_i^T\mathbf{Z}_i^{-1}\right\| \leq \min~\sum_{i\in \mathcal{N}}\lambda_{\mbox{\footnotesize max}}(\mathbf{Z}_i^{-1})=\max~\sum_{i\in \mathcal{N}}\lambda_{\mbox{\footnotesize min}}(\mathbf{Z}_i).
\end{eqnarray}
Equations (\ref{eq: 350}) and (\ref{eq: 251}) give (\ref{eq: 39-1}) with (D1) and (D3) and we can find the bound as
\begin{align}
\label{eq: 252}
\kappa_i\geq\|\mathbf{K}_i\|=\frac{1}{2}\|\mathbf{B}_i^T\mathbf{Z}_i^{-1}\|=\underline{\kappa_i}.
\end{align}
Derivation of (D2) and (D4) in (\ref{eq: 39-1}) and $\underline{\mu_i}$ is similar with $\mathbf{O}=\mathbf{0}$ and $\mathbf{M}=-\frac{1}{2}\hat{\mathbf{P}}^{-1}\mathbf{C}^T$.
\end{IEEEproof}

\emph{Corollary 1}: If the network is a partially ordered set (poset) then, decentralized control is possible if $\kappa_i\ge \underline{\kappa_i}$ and $\mu_i\ge \underline{\mu_i}$ where $\underline{\kappa_i}$ and $\underline{\mu_i}$ are given in (\ref{eq: 239}) with $\mathbf{H}=\mathbf{0}$ in (\ref{eq: 39-1}).

\begin{IEEEproof}
If the network is a poset, $\mathbf{H}$ is a block lower triangular matrix. Let us also limit $\mathbf{L}$ and $\mathbf{O}$ to have the same structure as $\mathbf{H}$. Therefore, the closed loop system matrices $\mathbf{A}+\mathbf{H}+\mathbf{B}(\mathbf{K}+\mathbf{L})$ and $\mathbf{A}+\mathbf{H}+(\mathbf{M}+\mathbf{O})\mathbf{C}$ in (\ref{eq: 3}) and (\ref{eq: 4}) are also block lower traingular. Hence, the closed loop eigenvalues are determined solely by matrices on the main diagonal, namely $\mathbf{A}_i+\mathbf{B}_i\mathbf{K}_i$ and $\mathbf{A}_i+\mathbf{M}_i\mathbf{C}_i$. Eigenvalues of these matrices can be arbitrarily placed in the open left half plane due to controllability of $(\mathbf{A}_i,\mathbf{B}_i)$ and observability of $(\mathbf{A}_i,\mathbf{C}_i)$. Thus, (\ref{eq: 7}) can be written as
\begin{align}
\left(\mathbf{A}+\mathbf{B}\mathbf{K}\right)^{T}\mathbf{P}+\mathbf{P}\left(\mathbf{A}+\mathbf{B}\mathbf{K}\right)+2\boldsymbol \beta \circ\mathbf{P}&\prec \mathbf{0}\nonumber\\
\left(\mathbf{A}+\mathbf{M}\mathbf{C}\right)^{T}\hat{\mathbf{P}}+\hat{\mathbf{P}}\left(\mathbf{A}+\mathbf{M}\mathbf{C}\right)+2\boldsymbol \beta \circ\hat{\mathbf{P}}&\prec \mathbf{0}\nonumber\\
\mathbf{P}=\mbox{diag}(\mathbf{P}_i) &\succ \mathbf{0}\nonumber\\
\hat{\mathbf{P}}=\mbox{diag}(\hat{\mathbf{P}}_i)&\succ  \mathbf{0}.\nonumber
\end{align}
Following along the same lines as the proof of Theorem 3, the bounds on decentralized controller gains are given in (\ref{eq: 239}) by setting $\mathbf{H}=\mathbf{0}$ in (\ref{eq: 39-1}) which is always feasible due to controllability of $(\mathbf{A}_i,\mathbf{B}_i)$ and observability of $(\mathbf{A}_i,\mathbf{C}_i)$.
\end{IEEEproof}



In general the convex mixed-binary program (\ref{eq: 39}) is NP-hard. In the worst case, one has to solve $2^{E}$ convex problems, carrying an exhaustive search on the binary variables, where $E$ is the total number of edges in the plant network. While a variety of exact methods for convex mixed-binary programs are available 
\cite{grossmann2002}, their computational complexity is prohibitive for large networks, or when the calculation is to be repeated periodically. Here, we propose a simple suboptimal relaxation-thresholding approach, as follows:
\begin{enumerate}
\item{Initialize $\alpha_{ij}\leftarrow 1$ for all $i,j\in \mathcal{N},i\neq j$.}
\item{Solve the LMI feasibility problem (C1) - (C8) in (\ref{eq: 39}) to find $\mathbf{Z}_i$,$\mathbf{W}_i$,$\mathbf{Y}_{ij}$,$\hat{\mathbf{P}}_i$,$\hat{\mathbf{W}}_i$ and $\hat{\mathbf{Y}}_{ij}$. If this problem is feasible, $\alpha^{\dagger}_{ij}\leftarrow \alpha_{ij}$, otherwise go to step 5, unless the problem is infeasible at the first iteration, in which case there is no solution and search is terminated here.}
\item{Solve \eqref{eq: 39} with (C9) relaxed to $0\le \alpha_{ij}\le 1$ to obtain solution $\alpha^{(r)}_{ij}$ satisfying (C1) and (C2), where $\mathbf{Z}_i$,$\mathbf{W}_i$,$\mathbf{Y}_{ij}$,$\hat{\mathbf{P}}_i$,$\hat{\mathbf{W}}_i$ and $\hat{\mathbf{Y}}_{ij}$ are those found in step 2.}
\item{If all $\alpha^{(r)}_{ij}=0$, set $\alpha^{\dagger}_{ij}\leftarrow 0$ and go to step 5. Otherwise, set $\alpha_{ij}$ corresponding to the smallest non-zero $\alpha^{(r)}_{ij}$ to zero and return to step 2.}
\item{Return $\alpha^{\dagger}_{ij}$.}
\end{enumerate}
Note that in the worst case, one has to solve $E$ convex problems, since it can be solved by using a linear search in a sorted set $\{\alpha_{ij}^{(r)}\}$.\

To further simplify the procedure we can substitute steps 4 and 5 above with
\begin{enumerate}
\setcounter{enumi}{3}
\item{Solve
\begin{eqnarray}
\label{eq:tau}
\begin{array}{rl}
\underset{k,l}{\text{maximize}}& \tau=\alpha^{(r)}_{kl}\\
\mbox{subject to}&  \alpha_{ij}=1_{\alpha^{(r)}_{ij}\geq \tau}, \mbox{(C1) and (C2)}.
\end{array}
\end{eqnarray}}
\item{Return $\alpha^{\dagger}_{ij}=1_{\alpha_{ij}^{(r)}>\tau^{\star}}$, where $\tau^\star$ is the solution of \eqref{eq:tau}.}
\end{enumerate}
We note that the maximum number of convex problems that should be solved in \eqref{eq:tau} is only $\log E$, since it can be solved by a binary search on $\tau$ in a sorted set $\{\alpha_{ij}^{(r)}\}$. Of course, this reduction in complexity is at the price of a more conservative solution.

\section{Numerical Example}\label{sec:example}
\begin{figure}[!t]
\centering
\includegraphics[scale=0.4]{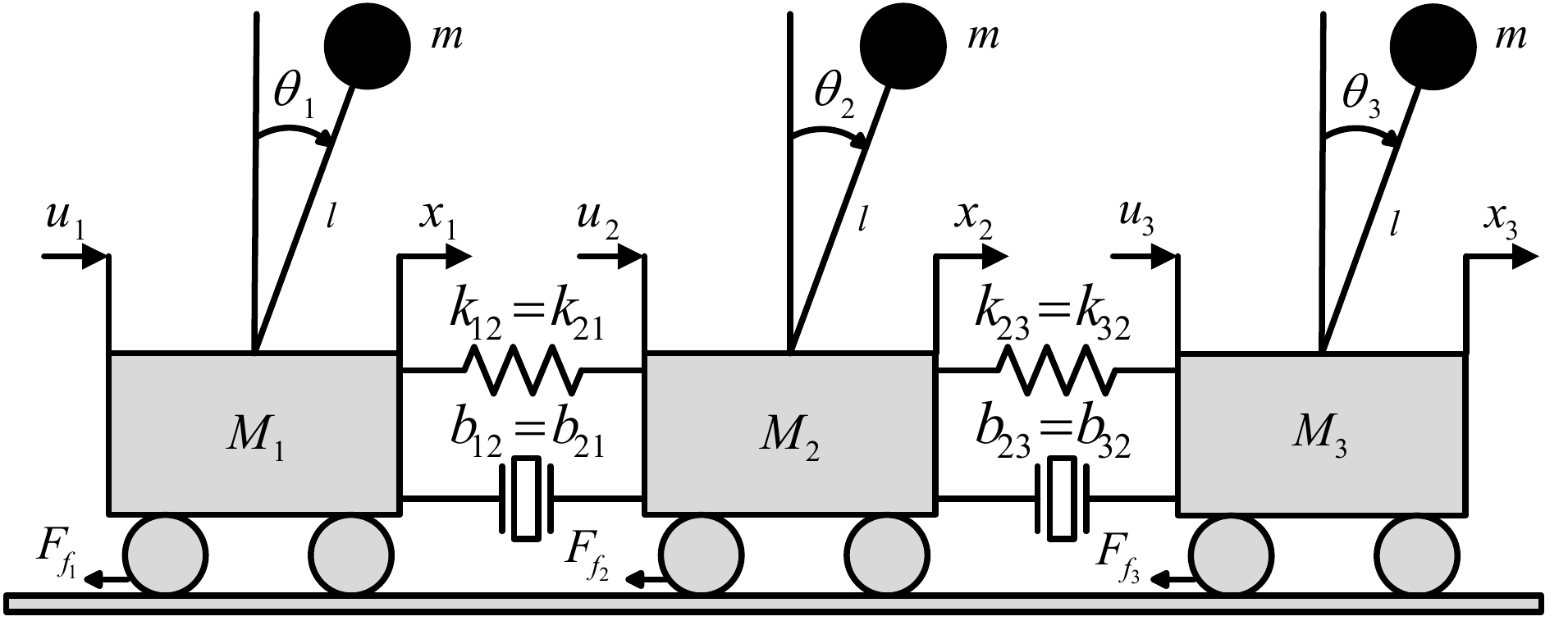}
\vspace{-0.1in}
\caption {Network of three coupled inverted pendulums}
\vspace{-0.2in}
\label{fig:NCS2}
\end{figure}
Consider the system shown in Fig. \ref{fig:NCS2}, where three inverted pendulums are mounted on coupled carts. The linearized equations of motions are \cite{ogata}
\begin{eqnarray}
\label{eq: 301}
M_il\ddot{\theta}_i &=& (M_i+m)g\theta_i+c_i\dot{x}_i+\sum_{j\in \mathcal{N}_i}\left[b_{ij}(\dot{x}_i-\dot{x}_j)+k_{ij}(x_i-x_j)\right]-u_i,\nonumber\\
M_i\ddot{x}_i &=& -c_i\dot{x}_i-\sum_{j\in \mathcal{N}_i}\left[b_{ij}(\dot{x}_i-\dot{x}_j)+k_{ij}(x_i-x_j)\right]-mg\theta_i+u_i.
\end{eqnarray}
where $c_i$, $b_{ij}=b_{ji}$, and $k_{ij}=k_{ji}$ are friction, damper and spring coefficients, respectively. Moreover, we have assumed that the moment of inertia of the pendulums is zero. Eq. (\ref{eq: 301}) can be written in the form of \eqref{eq: 1} by defining $\mathbf{x}_i=[x_{i,1}~x_{i,2}~x_{i,3}~x_{i,4}]^T=[\theta_i ~ \dot{\theta}_i ~ x_i ~ \dot{x}_i]^T$, $\mathbf{y}_i=[\theta_i ~ x_i]^T$,

\begin{align}
\label{eq: 302}
\mathbf{A}_i&=\left[
               \begin{array}{cccc}
                 0 & 1 & 0 & 0 \\
                 \frac{M_i+m}{M_il}g & 0 & \frac{k_i}{M_il} & \frac{c_i+b_i}{M_il} \\
                 0 & 0 & 0 & 1 \\
                 \frac{-m}{M_i}g & 0 & \frac{-k_i}{M_i} & \frac{-c_i-b_i}{M_i}
               \end{array}
             \right], \mathbf{H}_{ij}=\left[
          \begin{array}{cccc}
            0 & 0 & 0 & 0 \\
            0 & 0 & \frac{-k_{ij}}{M_il} & \frac{-b_{ij}}{M_il} \\
            0 & 0 & 0 & 0 \\
            0 & 0 & \frac{k_{ij}}{M_i} & \frac{b_{ij}}{M_i}
          \end{array}
        \right],\nonumber\\
\mathbf{B}_i&=\left[\begin{array}{cccc} 0 & \frac{-1}{M_il} & 0 & \frac{1}{M_i}\end{array}\right]^T, \mathbf{C}_i=\left[\begin{array}{cccc}1 & 0 & 0 & 0 \\0 & 0 & 1 & 0\end{array}\right],\nonumber
\end{align}
for $(i,j)=(1,2),(2,1),(2,3),(3,2)$, where $k_i=\sum_{j\in \mathcal{N}_i}k_{ij}$ and $b_i=\sum_{j\in \mathcal{N}_i}b_{ij}$. Since the subsystems are controllable and observable, we can apply Theorem 1 to design distributed observers and controllers that stabilize the entire network with small number of links in the control network. As design criteria, we assume each subsystem needs degree of stability $\beta_i=0.5$ and bounds on the norms of coupling gains are $\iota_{ij}=30$ and $\omega_{ij}=10$ and the numerical system parameters are $M_1=2$, $M_2=1$, $M_3=3$, $m=0.5$, $g=10$, $l=0.5$, $k_{12}=k_{21}=5$, $k_{23}=k_{32}=15$, $b_{12}=b_{21}=1$, $b_{23}=b_{32}=5$, $c_1=4$, $c_2=2$ and $c_3=1$. We will consider three different cases, where we progressively increase $\kappa_i$ and $\mu_i$. These parameters and the corresponding results are given in Table \ref{tab:res}, where $\alpha_{ij}^\dagger$ are found using the proposed simplified algorithm with linear search and $\alpha_{ij}^\star$ are found using an exhaustive search on the binary variables, followed by convex optimization of other variables. We can see that in this example the proposed relaxation-thresholding algorithm provides the optimal results.

We see that as the gain norm constraints are relaxed, the required control network becomes more sparse. In case 3 we see that if the bounds are relaxed beyond $\underline{\kappa_1}=54.1, \underline{\kappa_2}=273.2, \underline{\kappa_3}=152.1, \underline{\mu_1}=27.2, \underline{\mu_2}=29.2, \underline{\mu_3}=27.0$, decentralized control is possible, due to Theorem 3.

\begin{table*}[!t]
\caption{Numerical results for the three considered cases.}
\label{tab:res}
\vspace{-0.1in}
\centering
\begin{tabular}{|l|ccc|ccc|cccccc|cccccc|}\hline
& $\kappa_1$ & $\kappa_2$ & $\kappa_3$ & $\mu_1$ & $\mu_2$ & $\mu_3$ & \hspace{-0.05in}$\alpha_{12}^{\dag}$\hspace{-0.05in} & \hspace{-0.05in}$\alpha_{21}^{\dag}$\hspace{-0.05in} & \hspace{-0.05in}$\alpha_{13}^{\dag}$\hspace{-0.05in} & \hspace{-0.05in}$\alpha_{31}^{\dag}$\hspace{-0.05in} & \hspace{-0.05in}$\alpha_{23}^{\dag}$\hspace{-0.05in}& \hspace{-0.05in}$\alpha_{32}^{\dag}$\hspace{-0.05in} & \hspace{-0.05in}$\alpha_{12}^{\star}$\hspace{-0.05in} & \hspace{-0.05in}$\alpha_{21}^{\star}$\hspace{-0.05in} &\hspace{-0.05in}$\alpha_{13}^{\star}$\hspace{-0.05in} & \hspace{-0.05in}$\alpha_{31}^{\star}$\hspace{-0.05in} & \hspace{-0.05in}$\alpha_{23}^{\star}$\hspace{-0.05in}& \hspace{-0.05in}$\alpha_{32}^{\star}$\hspace{-0.05in}\\\hline\hline
1&96&106&211&27&26&28&1&1& 0 & 0 & 1&1&1&1& 0 & 0 &1&1\\\hline
2&135&121&232&27&28&29&0&0& 0 & 0 &1&1&0&0& 0 & 0 &1&1\\\hline
3&$\ge\underline{\kappa_1}$&$\ge\underline{\kappa_2}$&$\ge\underline{\kappa_3}$&$\ge\underline{\mu_1}$&$\ge\underline{\mu_2}$&$\ge\underline{\mu_3}$&0&0& 0 & 0 & 0 & 0 &0&0&0&0&0&0\\\hline
\end{tabular}
\end{table*}

\section{Concluding Remarks}
\label{sec:conclusion}
We have provided a design approach for distributed observer-based controllers that stabilize a given networked control system with an arbitrary directed topology. To measure states of each subsystem, we use the outputs of other nodes to improve stability of observer dynamics, in an approach dual to that of distributed controllers. Our design approach is based on a set of stability conditions obtained using the Lyapunov approach, and provides a sparse observer-controller network which guarantees global asymptotic stability. Moreover, we found conditions and bounds on norm of local controller gain matrices which allow complete decentralization. Due to some assumptions made to maintain tractability, the design includes some degree of conservatism. Thus, although the results provide us with significant insight into the problem of designing the sparsest controller-observer network, a gap still remains. Quantification or reduction of this gap will be quite valuable.

To avoid spending the entire margin in the stability criteria during the search for a sparse controller-observer network, we added a margin to the stability inequalities, as a free variable. Optimal distribution of this margin among the inequalities to make the network robust without significantly growing the size of the controller-observer network is, however, unknown. It is also interesting to understand the tradeoff between the stability margin and the sparsity of the observer-controller network.

We believe that the results presented in this paper provide a foothold for further progress towards understanding these interesting and important problems.

\bibliographystyle{IEEETran}

\end{document}